\documentclass[aop,citesort,MSNbibl,dvips]{arximspdf}

\doi{10.1214/10-AOP597}
\volume{39}
\issue{4}
\pubyear{2011}
\firstpage{1607}
\lastpage{1620}

\makeatletter

\newtheorem{prop}{Proposition}
\newproclaim{remark}{Remark}
\newproclaim{example}{Example}
\newtheorem{theorem}{Theorem}

\newproclaim{Assumption}{Assumption}
\newproclaim{Property}{Property}

\newcommand{\implies}{\Longrightarrow}
\newcommand{\mod}{\operatorname{mod}}

\newcommand{\E}{\mathrm{E}}

\newcommand{\Tr}{\operatorname{Tr}}

\newcommand{\X}{{X}}

\newcommand{\termj}{\mathbb{T}_{\mathbf j}}

\newcommand{\lldots}{\cdots}

\makeatother

\begin{document}
\begin{frontmatter}

\title{Convergence of joint moments for independent random patterned matrices}
\runtitle{Joint distribution}

\begin{aug}
\author[A]{\fnms{Arup} \snm{Bose}\thanksref{t1}\ead[label=e1]{bosearu@gmail.com}},
\author[A]{\fnms{Rajat Subhra} \snm{Hazra}\corref{}\ead[label=e2]{rajat\_r@isical.ac.in}} and
\author[B]{\fnms{Koushik} \snm{Saha}\ead[label=e3]{koushiksaha877@gmail.com}}
\runauthor{A. Bose, R. S. Hazra and K. Saha}
\affiliation{Indian Statistical Institute, Indian Statistical
Institute
and~Bidhannagar~Govt.~College}
\address[A]{A. Bose\\
R. S. Hazra\\
Statistics and Mathematics Unit\\
Indian Statistical Institute\\
203 B. T. Road, Kolkata 700108\\
India\\
\printead{e1}\\
\phantom{E-mail: }\printead*{e2}}
\address[B]{K. Saha\\
Department of Mathematics\\
Bidhannagar Govt. College\\
Salt Lake City, Kolkata 700064 \\
India\\
\printead{e3}}
\end{aug}

\thankstext{t1}{Supported by J. C. Bose National Fellowship,
Dept. of Science and Technology, Govt. of India. Part of the work was
done while visiting Dept. of Economics, University of Cincinnati.}

\received{\smonth{7} \syear{2010}}
\revised{\smonth{8} \syear{2010}}

%
\begin{abstract}
It is known that the joint limit distribution of independent Wigner
matrices satisfies a very special asymptotic independence, called
freeness. We study the joint convergence of a few other patterned
matrices, providing a framework to accommodate other joint laws. In
particular, the matricial limits of symmetric circulants and reverse
circulants satisfy, respectively, the classical independence and the
half independence. The matricial limits of Toeplitz and Hankel matrices
do not seem to submit to any easy or explicit independence/dependence
notions. Their limits are not independent, free or half independent.
\end{abstract}

%

\end{frontmatter}

\section{Introduction}
Wigner~\cite{wigner1958} showed how the semicircular law arises as the
limit of the empirical spectral distribution of a sequence of Wigner
matrices. See, for example,~\cite{Bai} and \cite
{andersonguionnetzeitouni} for such results and their variations. Then
researchers studied the joint convergence of independent Wigner
matrices and the limit is tied to the idea of free independence
developed by Voiculescu~\cite{voi1991}.

It appears that the study of joint distribution of random matrices has
been mostly concentrated on Wigner type matrices. In~\cite{ryan}, joint
limits of random matrices were studied and it was shown that there are
some circumstances where the limit may exist but
may not be free. As one of the Referees pointed out,~\cite{ryandebbah}
considers the joint distribution of Vandermonde matrices and diagonal
matrices and emphasizes the importance of studying joint distribution
of other patterned matrices.\vadjust{\goodbreak}

We study the joint convergence of $p$ independent
symmetric matrices with identical pattern. In particular, we show that
the tracial limit exists for any monomial when
the patterned matrix is any one of, Toeplitz, Hankel, symmetric
circulant or reverse circulant. The Wigner, symmetric circulant and
reverse circulant limits are, respectively, free semicircular, classical
independent normal and half independent Rayleigh. The Toeplitz and
Hankel limits are not free, independent or half independent.

In Section~\ref{ncp}, we discuss some preliminaries on noncommutative
probability spaces and recall the notions of independence, freeness and
half independence. In Section~\ref{pattern1}, we introduce the notion
of words and colored words
and state our main result. In Proposition~\ref{main:theo}, we show that
under a mild condition, if the marginal limit
exists then the joint limit exists and can be expressed in terms of the
marginals with the help of words and colored words. In Section~\ref
{examples}, we discuss some examples. Finally, in Section \ref
{proof}, we give a proof of Proposition~\ref{main:theo}.

\section{Noncommutative probability spaces and independence}\label{ncp}
A \textit{noncommutative probability space} is a pair $(\mathcal
{A},\phi
)$ where $\mathcal{A}$ is a unital complex algebra (with unity $1$) and
$\phi\dvtx\mathcal{A}\rightarrow\mathbb{C}$ is a linear functional
satisfying $\phi(1)=1$.

Two important examples of such spaces are the following:

\begin{longlist}[(1)]
\item[(1)] Let $(X,\mathcal{B}, \mu)$ be a probability space. Let $L(\mu
)=\bigcap_{1\leq p<\infty}L^p(\mu)$ be the algebra of random variables
with finite moments of all orders. Then $(L(\mu),\phi)$ becomes a
(commutative) probability space where
$\phi$
is the expectation operator, that is, integration with respect to $\mu$.

\item[(2)] Let $(X,\mathcal{B}, \mu)$ be a probability space and let
$\mathcal
{A}=\mathrm{Mat}_{n}(L(\mu))$ be the space of $n\times n$ complex random
matrices with elements from $L(\mu)$. Then $\phi$ equal to $\frac
{1}{n}\E_\mu[\Tr(\cdot)]$ or $\frac{1}{n}[\Tr(\cdot)]$ both yield
noncommutative probabilities.

For any noncommuting variables $x_1, \ldots, x_n$, let $\mathbb
{C}\langle x_1, x_2,\ldots,x_n \rangle$ be the unital algebra of all
complex polynomials in these variables. If $a_1, a_2, \ldots,a_n \in
\mathcal{A}$, then their joint distribution $\mu_{\{a_i\}}$ is defined
canonically by their mixed moments
$\mu_{\{a_i\}}(x_{i_1}\lldots x_{i_m}) =\phi(a_{i_1}\lldots a_{i_m})$.
That is,
\[
\mu_{\{a_i\}}(P)=\phi(P(\{a_i\})) \qquad\mbox{for } P \in\mathbb
{C}\langle
x_1, x_2,\ldots,x_n \rangle.
\]
\end{longlist}

\subsection*{Convergence in law} Let $(\mathcal{A}_n,\phi_n)$, $n
\geq1$, and
$(\mathcal{A},\phi)$ be noncommutative probability spaces and let $\{
a_i^n\}_{i \in J}$ be a sequence of subsets of $\mathcal{A}_n$ where
$J$ is any finite subset of
$\mathbb N$. Then we say that $\{a_i^n\}_{i \in J}$ \textit{converges in
law} to $\{a_i\}_{i \in J}\subset\mathcal{A}$ if for all complex
polynomials $P$,
\[
\lim_{n\rightarrow\infty}\mu_{\{a_i^n\}_{i\in J}}(P)=\mu_{\{a_i\}_{i
\in J}}(P).
\]
To verify convergence in the above definition, it is enough to verify
the convergence for all monomials $q=x_{i_1}\lldots x_{i_k}$, $k \geq1$.

\subsection*{Independence and free independence of algebras}
Suppose $\{{\mathcal{A}}_i\}_{i \in J}\subset{\mathcal{A}}$ are
unital subalgebras. They are said to be \textit{independent} if they
commute and $\phi(a_1\lldots a_n)=\phi(a_1)\lldots\phi(a_n)$ for all
$a_i \in\mathcal{A}_{k(i)}$ where ${i\neq j\implies k(i)\neq k(j)}$.

These subalgebras are called \textit{freely independent} or simply
\textit{free} if $\phi(a_j)=0$, $a_j \in{\mathcal{A}}_{i_j}$ and
$i_j\neq
i_{j+1}$ for all $j$ implies $\phi(a_1\lldots a_n)=0$.
The random variables (or elements of an algebra) $(a_1,a_2,\ldots)$
will be called independent (resp., free) if the subalgebras
generated by them are independent (resp., free).

\subsection*{Half independence of elements of an algebra}
Half independence
arises in classification results for easy quantum groups and some
quantum analogue of de Finetti's theorem. We closely follow the
developments in~\cite{banica2010}. We shall see later how this half
independence arises in the context of convergence of random matrices.
To describe this, we need the concepts of half commuting elements and
symmetric monomials.

\subsection*{Half commuting elements}
Let $\{a_i\}_{i \in J} \subset
\mathcal{A}$. We say that they
\textit{half commute} if $ {a_ia_ja_k=a_ka_ja_i}$, for all $ i, j, k
\in J$.
Observe that if $\{a_i\}_{i \in J}$ half commute then $a_i^2$ commutes
with $a_j$ and $a_j^2$ for all $i,j \in J$.

\subsection*{Symmetric monomials}
Suppose $\{a_i\}_{i \in J} \subset\mathcal{A}$. For any $k \geq1$,
and any $\{i_j\} \subset J$, let $a=a_{i_1}a_{i_2}\lldots a_{i_k}$ be an
element of $\mathcal{A}$. For any $ i \in J$, let $E_i(a)$ and $O_i(a)$
be, respectively, the number of times $a_i$ has occurred in the even
positions and in the odd positions in $a$.
The monomial $a$ is said to be \textit{symmetric} (with respect to $\{
a_i\}
_{i \in J}$) if $E_i(a)=O_i(a)$ for all $i \in J$. Else it is said to
be nonsymmetric.

\subsection*{Half independent elements}
Let $\{a_i\}_{i \in J}$ in $(\mathcal{A},\phi)$ be half commuting.
They are said to be half independent if (i) $\{a_i^2\}_{ i \in J}$ are
independent and (ii) whenever $a$ is nonsymmetric with respect to $\{
a_i\}_{i \in J}$, we have $\phi(a)=0$.
\begin{remark} The above definition is equivalent to that given
in~\cite{banica2010}, although there is no notion of symmetric
monomials there. As pointed out in~\cite{speicher1997}, the concept of
half independence does not extend to subalgebras.
\end{remark}
\begin{example}\label{exampl1} This
is Example 2.4 of~\cite{banica2010}.
Let $(\Omega,\mathcal{B},\mu)$ be a probability space and let $\{
\eta
_i\}$ be a family of independent complex Gaussian random variables.
Define $a_i \in(M_2(L(\mu)),\E[\Tr(\cdot)])$ by
\[
a_i= \left[
\matrix{
0 & \eta_i \cr
\bar\eta_i & 0
}
\right].
\]
These $\{a_i\}$
are half independent.
\end{example}
\begin{remark}
Let $X,Y$ and $Z$ be three self-adjoint elements of
$(\mathcal{A},\phi)$ such that $\phi(X)=\phi(Y)=\phi(Z)=0$ but
$\phi
(X^2),\phi(Y^2),\phi(Z^2)\neq0$. Following Remark 5.3.2 of \cite
{andersonguionnetzeitouni}:\vadjust{\goodbreak}

\begin{longlist}
\item If $X,Y$ commute and are independent, then $ \phi(XY)=0$ and\break $\phi
(XYXY)=\phi(X^2)\phi(Y^2)\neq0$.

\item
If $X$, $Y$ and $Z$ are half independent,
then $ {\phi(XYZZXY)=0.}$
This happens since $X$ appears two times in odd positions but zero
times in even positions.

\item
If $X,Y$, $Z$ are free, then $\phi(XYZXYZ)=0$. If they
are half independent then $\phi(XYZXYZ)=\phi(X^2)\phi(Y^2)\phi
(Z^2)\neq
0$.
\end{longlist}
\end{remark}

\section{Joint convergence of patterned matrices}\label{pattern1}
\subsection{Some preliminaries and the main result}
It is well known that the joint limit of independent Wigner matrices
yields the freely independent semicircular law.
At the same time, there are a host of results on the limiting spectrum
of other matrices. Important examples include the sample variance--covariance,
the Toeplitz and the Hankel matrices. The joint convergence
of several sample variance--covariance matrices has also been
investigated in the literature. However, joint convergence does not
seem to have been addressed in any generality.
In particular, it is not clear what other notions of independence are
possible when we consider the joint limit of (independent) matrices.

Patterned matrices offer a general framework for which this question
may be worth investigating.
A general pattern matrix may be defined through a \textit{link function}.
Let $d$ be a positive integer. Let $\mathbb Z$ be the set of all
integers and let $\mathbb N$ be the set of all natural numbers. Let
$L_n\dvtx\{1, 2, \ldots, n\}^2 \to\mathbb Z^d, n \geq1$, be a
sequence of
functions such that $L_{n+1}(i,j)=L_n
(i,j)$ whenever $ 1\leq i, j \leq n$. We shall write $L_n=L$ and call
it the \textit{link} function and
by abuse of
notation we write $\mathbb N^2$ as the common domain of $\{ L_n \}$.
For our examples later, the value of $d$ is either
$1$ or $2$.

A typical patterned matrix is then of the form $A_n=(x(L(i,j)))$ where
$\{ x (i); i \geq0\}$ or $\{ x(i, j); i, j\geq1 \}$ is a sequence of
variables. In what follows, we only consider real symmetric matrices.
Here are some well-known matrices and their link functions:

\begin{longlist}
\item
Wigner matrix $W_n$. $ L\dvtx\mathbb N^2 \rightarrow\mathbb Z^2$ where
$L(i,j) = (\min(i,j) , \max(i,j))$.
\item
Symmetric Toeplitz matrix $T_n$. $ L\dvtx\mathbb N^2 \rightarrow
\mathbb Z $ where
$L(i,j) =|i-j|$.
\item
Symmetric Hankel matrix $H_n$. $L\dvtx\mathbb N^2 \rightarrow
\mathbb
Z$ where
$L(i,j) =i+j$.
\item
Reverse circulant matrix $\mathit{RC}_n$. $L\dvtx\mathbb N^2 \rightarrow
\mathbb Z $ where
$L(i,j) =(i+j )\mod n$.
\item
Symmetric circulant matrix $\mathit{SC}_n$. $L\dvtx\mathbb N^2 \rightarrow
\mathbb Z $ where
$L(i,j) =n/2-|n/2-|i-j||$.
\end{longlist}

In general, we assume that the link function $L$ satisfies Property
\ref{properB}.
\renewcommand{\theProperty}{B}
\begin{Property}\label{properB}
$\Delta(L) = \sup_n \sup_{ t \in\mathbb Z^d } \sup
_{1 \le k \le n} \# \{
l\dvtx1 \le l \le n, L(k,l) = \break t\} < \infty$.

In particular, $\Delta(L)=2$ for $T_n,\mathit{SC}_n$, and $\Delta(L)=1$ for
$W_n,H_n$ and $\mathit{RC}_n$.
\end{Property}

Now let $(\Omega,\mathcal{B},\mu)$ be a probability space and let
$X_{i,n}\dvtx\Omega\rightarrow M_n$ for $1\leq i \leq p$ be symmetric
patterned random matrices of order $n$.
We shall refer to the $p$ indices as $p$ distinct \textit{colors}. The
$(j,k)$th entry of the matrix $X_{i,n}$ will be denoted by
$X_{i,n}(L(j,k))$.
\renewcommand{\theAssumption}{I}
\begin{Assumption}\label{AssumpI}
Let the input sequence of each matrix in the
collection $\{X_{i,n}\}_{1\leq i \leq p}$ be independent with mean zero
and variance $1$ and assume they are also independent across $i$.
Suppose the matrices have a common link function $L$ which satisfies
Property~\ref{properB} and
\[
\sup_{n \in\mathbb{N}}\sup_{1\leq i\leq p}\sup_{1\leq m\leq l\leq
n}\E
[|X_{i,n}(L(m,l))|^k]\leq c_k<\infty.
\]
Suppressing\vspace*{1pt} the dependence on $n$, we shall simply write $X_i$ for $X_{i,n}$.
We view $\{\frac{1}{\sqrt n}X_{i}\}_{1\leq i \leq p}$ as elements of
$(\mathcal{A}_n=\mathrm{Mat}_n(L(\mu)),\phi_n)$ where $\phi_n=n^{-1}\E[\Tr]$.
Denote the joint distribution of $\{\frac{1}{\sqrt n}X_{i}\}_{1\leq i
\leq p}$ by $\widehat{\mu}_n$.
Then $\{\frac{1}{\sqrt n}X_{i}\}_{1\leq i \leq p}$ converges in law if
%
%
\begin{eqnarray}\quad
\widehat{\mu}_n(q)&=&\phi_n(q)\nonumber\\[-8pt]\\[-8pt]
&=&\frac
{1}{n^{1+k/2}}\E[\Tr
(X_{i_1}\lldots X_{i_k})]\nonumber\\
\label{equation2}
&=&\frac{1}{n^{1+k/2}}\sum_{j_1,\ldots,j_k}\E
[X_{i_1}(L(j_1,j_2))X_{i_2}(L(j_2,j_3))\lldots X_{i_k}(L(j_k,j_1))]
\end{eqnarray}
converges for all monomials $q$ of the form $q(\{X_i\}_{1\leq i \leq
p})=X_{i_1}\lldots X_{i_k}$.
\end{Assumption}

To upgrade to almost sure convergence, we also define
\begin{eqnarray*}
\widetilde{\mu}_n(q)&=&\frac{1}{n^{1+k/2}}\Tr[X_{i_1}\lldots
X_{i_k}]\\
&=&\frac
{1}{n^{1+k/2}}\sum_{j_1,\ldots
,j_k}X_{i_1}(L(j_1,j_2))X_{i_2}(L(j_2,j_3))\lldots X_{i_k}(L(j_k,j_1)).
\end{eqnarray*}
All developments below are with respect to one fixed monomial $q$ at a
time.

\subsubsection*{Circuit} Any $\pi\dvtx\{0,1,2,\ldots,h\}
\rightarrow
\{1,2,\ldots,n\} $ with $\pi(0) =
\pi(h)$
is
a \textit{circuit} of
\textit{length} $l(\pi):=h$. The dependence of a circuit on $h$ and $n$
will be suppressed.
A~typical element in (\ref{equation2}) can be now written as
%
%
\begin{equation}\label{equation3}
\E\Biggl[\prod_{j=1}^k X_{i_j}\bigl(L\bigl(\pi(j-1),\pi(j)\bigr)\bigr)\Biggr].
\end{equation}
If all $L$-values $ L(\pi(j-1),\pi(j))$ are repeated more than once
in (\ref{equation3}), then
the circuit is \textit{matched}. If $L$ values are repeated exactly
twice, then it is called \textit{pair matched}.
If the $L$ values are repeated within the same color, then it is
\textit{color matched.}\vadjust{\goodbreak}

For $q=X_{i_1}X_{i_2}\lldots X_{i_k}$, let for convenience, the
corresponding sequence of colors be denoted by $\{c_1,c_2,\ldots,c_k\}
$. Also let
$H=\{\pi\dvtx\pi\mbox{ is a color matched}\break \mbox{circuit}\}$. Define
an equivalence relation on $H$ by defining $\pi_1\sim_C \pi_2$ if and
only if, $c_i=c_j$ and
\begin{eqnarray*}
&&X_{c_i}\bigl(L\bigl(\pi_1(i-1),\pi_1(i)\bigr)\bigr)=\X_{c_j}\bigl(L\bigl(\pi_1(j-1),\pi
_1(j)\bigr)\bigr)\\
&&\quad\Longleftrightarrow\quad X_{c_i}\bigl(L\bigl(\pi_2(i-1),\pi_2(i)\bigr)\bigr)=\X
_{c_j}\bigl(L\bigl(\pi_2(j-1),\pi_2(j)\bigr)\bigr).
\end{eqnarray*}

\subsubsection*{Colored words}
An equivalence class induces a partition of $\{1, 2,\ldots,k\}$ and
each block of the partition is associated with a color. Any such class
can be expressed as a (colored) word $w$ where letters appear in
alphabetic order of their first occurrence and with a subscript to
distinguish the color.
For example, the partition $(\{\{1,3\},1\}, \{\{2,4\},2\}, \{\{5,7\}
,1\},\{\{6,8\},3\})$ is identified with the word
$a_1b_2a_1b_2c_1d_3c_1d_3$. A typical position in a colored word would
be referred to as $w_{c_i}[i]$.

Let the class of all (colored) circuits corresponding to a color
matched word $w$ and the class of all pair matched colored words be
denoted, respectively,
by
%
%
\begin{eqnarray}\quad
\Pi_{C_q}(w)&=&\bigl\{\pi\dvtx
w_{c_i}[i]=w_{c_j}[j]\Longleftrightarrow X_{c_i}\bigl(L\bigl(\pi(i-1),\pi
(i)\bigr)\bigr)\nonumber\\[-4pt]\\[-12pt]
&&\hspace*{98.3pt}=X_{c_j}\bigl(L\bigl(\pi(j-1),\pi(j)\bigr)\bigr)\bigr\},\nonumber\\
\mathit{CW}_k(2)&=&\{\mbox{all paired matched (within same color) words }w
\mbox
{ of length } k
\} \nonumber\\[-4pt]\\[-12pt]
&&\eqntext{\mbox{($k$ is even)}.}
\end{eqnarray}
All the above notions have the corresponding noncolored versions. For
instance, if we drop the colors from a colored word, then we obtain a
\textit{noncolored} word. For any monomial $q$, dropping the color
amounts to dealing with only one matrix or in other words with the
marginal distribution.

Let $w[i]$ denote the $i$th entry of a noncolored word
$w$. The equivalence class corresponding to $w$ and the set of pair
matched noncolored words will be denoted, respectively, by
%
%
\begin{eqnarray} \Pi(w) &=& \bigl\{\pi\dvtx w[i] = w[j] \Leftrightarrow
L\bigl(\pi
(i-1),\pi(i)\bigr) =
L\bigl(\pi(j-1),\pi(j)\bigr)\bigr\},\\
W_k(2)& = &\{\mbox{all paired matched words } w \mbox{ of length }
k\}\qquad
\mbox{($k$ is even)}.
\end{eqnarray}
For any word $w \in \mathit{CW}_k(2)$, consider the noncolored word $w'$
obtained by dropping the color.
Then $w' \in W_k(2)$. Since we are dealing with one fixed monomial at
a time, this yields a bijective mapping say
%
%
\begin{equation}\label{one-one mapping}\psi_q\dvtx \mathit{CW}_k(2)
\rightarrow W_k(2).
\end{equation}

For any $w\in \mathit{CW}_k(2)$, define
\[
{p_{C_q}(w)=\lim_{n\rightarrow\infty}\frac{1}{n^{k/2+1}}|\Pi_{C_q}(w)|}\qquad
\mbox{if the limit exists}.
\]
\begin{prop}\label{main:theo}
Let $\{X_i\}_{1\leq i \leq p}$ be
patterned matrices
satisfying Assumption~\ref{AssumpI}. Fix any monomial
$q=X_{i_1}X_{i_2}\lldots
X_{i_{k}}$. Assume that, whenever~$k$ is even,
%
%
\begin{equation}\label{wordexistence}p(w)=\lim_{n\rightarrow\infty
}\frac{1}{n^{k/2+1}}|\Pi(w)|\qquad
\mbox{exists for all } w \in W_k(2).
\end{equation}

\textup{(a)} Then $p_{C_q}(w)=p(\psi_q(w))$ and for any $k$,
\[
\lim_{n\rightarrow\infty}\widehat{\mu}_n(q)=\sum_{w\in
\mathit{CW}_k(2)}p_{C_q}(w)=\alpha(x_{i_1}\lldots x_{i_{k}})\qquad
(\mbox{say})
\]
with
\begin{eqnarray*}
|\alpha(x_{i_1}\lldots x_{i_{k}})|&\leq&\frac{k!\Delta
(L)^{k/2}}{(k/2)!2^{k/2}}\qquad \mbox{if $k$ is even and each
color}\\[-4pt]
&&\hphantom{\frac{k!\Delta
(L)^{k/2}}{(k/2)!2^{k/2}}}\qquad\mbox{appear even number of times}\\
&=&0 \qquad\mbox{if $k$ is odd or a color appears odd number of times}.
\end{eqnarray*}

\textup{(b)} $ \E[|\widetilde{\mu}_n(q)-\widehat{\mu}_n(q)|^4]=O(n^{-2})$
and hence $ {\lim_{n\rightarrow\infty} \widetilde{\mu}_n(q)=\alpha
(x_{i_1}\lldots x_{i_{k}})}$ almost surely.
\end{prop}
\begin{remark}\label{existenceremark} It is\vspace*{1pt} known from~\cite{bosesen}
that (\ref{wordexistence}) holds true for Wigner, Toeplitz, Hankel,
reverse circulant and symmetric circulant matrices. The quantity
$\frac{k!}{(k/2)!2^{k/2}}$ above is the total number of noncolored
pair matched words of length $k$ ($k$ even). Often not all pair matched
words contribute to the limit and in such cases, this bound can be improved.

Consider the polynomial algebra $\mathbb{C}\langle a_1, a_2,\ldots,
a_p\rangle$ in noncommutative indeterminates $\{a_i\}_{1\leq i\leq p}$
and define a linear functional $\phi$ on it by
\[
\phi(a_{i_1}\lldots a_{i_k})=\lim_{n\rightarrow\infty}\widehat{\mu}
_n(X_{i_1}\lldots X_{i_k}).
\]
Then Proposition~\ref{main:theo} implies we have convergence in law of
$(\{\frac{1}{\sqrt n}X_{i}\}_{1\leq i \leq p}, \phi_n)$ to $(\{a_i\}
_{1\leq i\leq p}, \phi)$ where $\phi_n$ equals
$\frac{1}{n}\E[\Tr(\cdot)]$ or $\frac{1}{n}[\Tr(\cdot)]$. In the latter
case, the convergence is almost sure.
\end{remark}
\begin{remark}
Proposition~\ref{main:theo} shows that the joint moments of pattern matrices
can be expressed as functions of pair matched words or, in other
words, pair partitions.
\cite{bazekospeicher} consider bosonic, fermionic and $q$-Brownian
motions and show that the joint distribution of certain operators (in
some appropriate sense) can be expressed as functions on the set of
pair partitions. It would be interesting to investigate if there are
connections between the two types of models.
\end{remark}

\subsection{Some examples}\label{examples}
From the above result,
$\lim\widehat{\mu}_n(q)=0$ when $k$ is odd or when there is a color
which appears an odd number of times in the monomial $q$.
Henceforth, we thus assume
that the order of the monomial is even and each color appears an even
number of times.\vadjust{\goodbreak}
\begin{example}[(Wigner matrices)]\label{exampl2} Joint convergence of the
Wigner matrices was first studied in~\cite{voi1991} and later many
authors extended it. For details of the classical proof and further
extensions, we refer the readers to~\cite{andersonguionnetzeitouni}.
Here we give a quick partial proof essentially translating the concept
of noncrossing partitions that is used in the standard proof into words.
\end{example}

\subsubsection*{Colored Catalan words} Fix $k \geq2$. If for a $w \in
\mathit{CW}_k(2)$,
sequentially~deleting all double letters of the same color leads to the
empty word then we call $w$ a \textit{colored Catalan word}.
For example, the monomial $X_1X_2X_2X_1X_1X_1$ has exactly two colored
Catalan words $a_1b_2b_2a_1c_1c_1$ and $a_1b_2b_2c_1c_1a_1$. A colored
Catalan word associated with $X_1X_2X_2X_1X_1X_1X_2X_2$ is
$a_1b_2b_2a_1c_1c_1d_2d_2$ which is not even a valid colored word for
the monomial $X_1X_1X_1X_1X_2X_2X_2X_2$.

Let $\{W_{i}\}_{1\leq i \leq p}$ be an independent sequence of $n
\times n$ Wigner matrices satisfying Assumption~\ref{AssumpI}. Then from
Proposition~\ref{main:theo} and Remark~\ref{existenceremark},
$\{n^{-1/2}W_i\}_{1\leq i \leq p}$ converges in law to $\{a_i\}_{1\leq
i \leq p}$. We show that
$\{a_i\}_{1\leq i \leq p}$ are free and the marginals are distributed
according to the semicircular law.

From Table 1 of~\cite{bosesen}, for noncolored words, $p(w)$ equals $1$
if $w$ is a Catalan word and otherwise $p(w)=0$. As a consequence, the
marginals are semicircular.
Now fix any monomial
$q=x_{i_1}x_{i_2}\lldots x_{i_{2k}}$ where each color appears an even
number of times in the monomial.

Let $w$ be a colored Catalan word. It remains Catalan when we ignore
the colors. Hence from above,
$p_{C_q}(w)=p(\psi_q(w))=1$. Likewise, if $w$ is not colored Catalan
then the word $\psi_q(w)$ cannot be Catalan and hence
$p_{C_q}(w)=p(\psi
_q(w))=0$.
Hence, if $\mathrm{CAT}_q$ denotes the set of colored Catalan words
corresponding to a monomial $q$ then from the above discussion,
\[
\lim_{n\rightarrow\infty}\widehat{\mu}_n(q)= |\mathrm{CAT}_q|.
\]
Any double letter corresponds to a pair partition (within the same
color) by the equivalence relation $\sim_C$. It is known that the
number of Catalan words of length $2k$ is same as the number of
noncrossing pair partitions [denoted by $\mathit{NC}_{2k}(2)$] of length $2k$.
See~\cite{bosesen} and Chapter 1 of~\cite{andersonguionnetzeitouni}
for proofs.

Since the elements of the same pair partition must belong to the same
color, we have
\[
|\mathrm{CAT}_q|= \sum_{\pi\in \mathit{NC}_{2k}(2)}\prod_{(j,j') \in\pi}\mathbb
{I}_{c_j=c_{j'}}.
\]
This is precisely the free joint semicircular law corresponding to $q$
(see Theorem~5.4.2 of~\cite{andersonguionnetzeitouni}).

Incidentally, since the number of noncolored Catalan words of length
$2k$ is $\frac{2k!}{k!(k+1)!}$,
we have $|\alpha(x_{i_1}\lldots x_{i_{2k}})|\leq\frac{2k!}{k!(k+1)!}$.
Corollary 5.2.16 of~\cite{andersonguionnetzeitouni} can hence be
applied to claim the existence of a $C^*$-probability space with a
state $\phi$ and free semicircular random variables $\{a_i\}$ in it.
\begin{example}[(Symmetric circulants)]\label{exampl3}
The case of symmetric circulant is rather easy.
These matrices are commutative and so the limit is also commutative.

Let\vspace*{1pt} $\{\mathit{SC}_{i}\}_{1\leq i \leq p}$ be an independent sequence of $n
\times n$ symmetric circulant matrices satisfying Assumption \ref
{AssumpI}. Then
$ {\{n^{-1/2}\mathit{SC}_i\}_{1\leq i \leq p}}$ converges in law to $\{a_i\}
_{1\leq i\leq p}$ which are independent and the marginals are
distributed according to the standard Gaussian law.
To see this, first recall that the total number of (noncolored) pair
matched words of length $2k$ equals
$\frac{2k!}{k!2^k}= C_k \mbox{ (say)}$. Further
(see~\cite{bosesen}), for \textit{any} pair matched word $w\in W_{2k}(2)$,
\[
p(w)=\lim_n \frac{1}{n^{1+k}} |\Pi(w)| = 1.
\]
Now consider an order $2k$ monomial where each color appears an even
number of times.
Hence, from Proposition~\ref{main:theo}, for any fixed monomial $q$,
\[
\lim_{n\rightarrow\infty}\widehat{\mu}_n(q)=\sum_{w \in
\mathit{CW}_{2k}(2)}p_{C_q}(w)=|\mathit{CW}_{2k}(2)|.
\]
Let $l$ be the total number of distinct colors (distinct matrices) in
the monomial
$q=x_{i_1}x_{i_2}\lldots x_{i_{2k}}$. Let $2\times n_i$ be the number
of matrices of the $i$th color.
Then the set of all pair matched colored words of length $2k$ is
obtained by forming pair matched subwords of color $i$ of lengths
$2n_i$, $1 \leq i \leq l$. Hence,
%
%
\begin{equation}\label{moment:symmetriccirculant}\phi(a_{i_1}\lldots
a_{i_{2k}})= \prod_{i=1}^lC_{n_i}.
\end{equation}
Thus if $\{a_1,\ldots, a_p\}$ denotes i.i.d. standard normal random
variables, then the above is the mixed moment
$\E[\prod_{i=1}^l a_i^{2n_i}]$.
\end{example}
\begin{example}[(Reverse circulant)]\label{exampl4}
It can be easily observed using the link function that the reverse
circulant matrices are half commuting. This motivates the next theorem.
\end{example}
\begin{theorem}\label{reversecirculant:halfindependent}
Let $\{\mathit{RC}_{i}\}_{1\leq i \leq p}$ be an independent sequence of $n
\times n$ reverse circulant matrices satisfying Assumption \ref
{AssumpI}. Then
$ {\{n^{-1/2}\mathit{RC}_i\}_{1\leq i \leq p}}$ converges in law to half
independent $\{a_i\}_{1\leq i\leq p}\in(M_2(L(\mu)),\E[\Tr(\cdot)])$
where $a_i=\bigl[{0\atop\bar\eta_i}\enskip{\eta_i\atop0}\bigr]$ and $\eta
_i$ are i.i.d. complex Gaussian.
\end{theorem}

To prove the result, we need the following notion.

\subsubsection*{Colored symmetric words} Fix $k\geq2$. A word $w \in
\mathit{CW}_k(2)$ is
called \textit{colored symmetric} if each letter occurs once each in an
odd and an even position \textit{within the same color}. Clearly, every
colored Catalan word is a colored symmetric word.
\begin{pf*}{Proof of Theorem~\ref{reversecirculant:halfindependent}}
Consider a monomial $q$ of length $2k$ where each color appears an even
number of times.
From the single matrix case, it follows that $p(w)=0$ if $w$ is not a
symmetric word (see Table 1 of~\cite{bosesen}). If $w$ is not a colored
symmetric word, then $\psi_q(w)$ is not a symmetric word and hence for
such $w$, $p_{C_q}(w)=p(\psi_q(w))=0$. Hence, we may restrict to
colored symmetric words and then we have by Proposition \ref
{main:theo}(a) that
\[
\lim_{n\rightarrow\infty}\widehat{\mu}_n(q)=|\mathit{CS}_q(w)|,
\]
where $\mathit{CS}_q(w)$ is the collection of all colored symmetric words of
length $2k$.

The number of symmetric words of length $2k$ is $k!$. Let, as before,
$l$ be the number of distinct colors in the monomial and $2 n_i$ be the
number of matrices of the $i$th color. All symmetric words are
obtained by arranging the $2n_i$ letters of the $i$th color in a
symmetric way for $i=1, 2, \ldots, l$.

It is then easy to see that these arguments imply that
%
%
\begin{equation}\label{symfactor}|\mathit{CS}_q(w)|= n_1!\times n_2!\times
\cdots\times n_l!.
\end{equation}
First, observe that if the monomial $a_{i_1}a_{i_2}\lldots a_{i_k}\in
(M_2(L(\mu)),\E[\Tr(\cdot)])$ is nonsymmetric, then
\[
\E(\Tr[a_{i_1}a_{i_2}\lldots a_{i_k}])=0.
\]
If instead $q(\{a_i\})=a_{i_1}a_{i_2}\lldots a_{i_{2k}}$ is symmetric,
then we have by half independence (Example~\ref{exampl1}),
\[
\E(\Tr[a_{i_1}\lldots a_{i_{2k}}])=n_1!\times n_2!\times\cdots\times
n_l!=\lim_{n\rightarrow\infty}\widehat{\mu}_n(q).
\]
So it follows from (\ref{symfactor}) that the joint limit is
asymptotically half independent. Incidentally the moments $\{k!, k\geq
1\}$ are the $(2k)$th moments of the
symmetrized Rayleigh distribution.
\end{pf*}
\begin{example}[(Toeplitz and Hankel)]\label{exampl5}
Consider first the Toeplitz matrix. Since $p(w)$ exists, from
Proposition~\ref{main:theo}, we have the joint convergence for Toeplitz
matrices.
For any fixed monomial $q$, let $\mathit{SNC}_q$ be the colored symmetric words
which are not Catalan. Then we obtain the following:
\begin{eqnarray*}
\phi(a_{i_1}\lldots a_{i_k})&=&\sum_{w \in \mathit{CW}_k(2)}p_{C_q}(w)\\[-2pt]
&=&\sum_{w \in \mathrm{CAT}_q}p_{C_q}(w)+\sum_{w \in \mathit{SNC}_q}p_{C_q}(w)\\[-2pt]
&&{}+\sum
_{\mathrm{other}\ \mathrm{pair}\ \mathrm{matched}\ \mathrm
{colored}\
\mathrm{words}}p_{C_q}(w)\\[-2pt]
&=&|\mathrm{CAT}_q|+\sum_{w \in \mathit{SCN}_q}p_{C_q}(w)\\[-2pt]
&&{}+\sum_{\mathrm{other}\ \mathrm
{pair}\ \mathrm{matched}\ \mathrm{colored}\ \mathrm{words}} p_{C_q}(w).
\end{eqnarray*}
Consider $q(X_1,X_2,X_3)=X_1X_2X_3X_1X_2X_3$ where $X_1,X_2$ and $X_3$
are scaled independent Toeplitz matrices.
From
Table 4 of~\cite{bosesen},
\[
p_{C_q}(a_1b_2c_3a_1b_2c_3)=p(\psi
_q(a_1b_2c_3a_1b_2c_3))=p(abcabc)=\tfrac{1}{2}.
\]
For this monomial, the only pair matched colored word possible is
$a_1b_2c_3a_1b_2c_3$ and hence
$\phi(a_1a_2a_3a_1a_2a_3)=\frac{1}{2}\neq0$. Thus, the limit is not free.

Now let $q(X_1,X_2,X_3)=X_1X_2X_3X_2X_3X_1$.
Then the only pair of matched colored word is $a_1b_2c_3b_2c_3a_1$ and
$\phi(a_1a_2a_3a_2a_3a_1)=p_{C_q}(a_1b_2c_3b_2c_3a_1)=p(abcbca)=\frac{2}{3}$.
On the other hand we have already seen that $\phi
(a_1a_2a_3a_1a_2a_3)=\frac{1}{2}$.
Since the two contributions are not equal, the Toeplitz limit is not
independent.

If they had been half independent, then $\phi
(a_1a_2a_3a_1a_2a_3)=\phi(a_1^2)\phi(a_2^2)\times\phi(a_3^2)=1$, but that is
not the case.
Thus,
the Toeplitz limit is not free, independent or half independent.

For Hankel matrices, the colored nonsymmetric words \textit{do not}
contribute to the limit. So for any fixed monomial $q$ we have
\[
\phi(a_{i_1}\lldots a_{i_k})= |\mathrm{CAT}_q|+\sum_{w \in \mathit{SNC}_q}p_{C_q}(w).
\]
That the Hankel limit is also not free, half independent or independent
can be checked along the above lines by considering appropriate
monomials and their contributions. It is interesting to note that
Hankel matrices do not half commute and that is why even though the
limits vanish on nonsymmetric words they are not half independent.
\end{example}

\subsection{\texorpdfstring{Proof of Proposition \protect\ref{main:theo}}{Proof of Proposition 1}}\label{proof}

(a) Fix a monomial $q=q(\{X_i\}_{1\leq i\leq p})=X_{i_1}\cdots X_{i_k}$.
Since $\psi_q$ is a bijection,
\[
\Pi_{C_q}(w)=\Pi(\psi_q(w)) \qquad\mbox{for } w \in \mathit{CW}_k(2).
\]
Hence using (\ref{wordexistence}),
\[
\lim_{n\rightarrow\infty}\frac{1}{n^{k/2+1}}|\Pi_{C_q}(w)|=\lim
_{n\rightarrow\infty}\frac{1}{n^{k/2+1}}|\Pi(\psi_q(w))|
=p(\psi_q(w))= p_{C_q}(w).
\]
For simplicity, denote
\[
\mathbb{T}_{\mathbf j}=\E
[X_{i_1}(L(j_1,j_2))X_{i_2}(L(j_2,j_3))\cdots
X_{i_k}(L(j_k,j_1))] \qquad\mbox{for } \mathbf j=(j_1,\ldots,j_k) .
\]
Then
%
%
\begin{equation}\label{tj}\widehat{\mu}_n(q)=\frac
{1}{n^{k/2+1}}\sum
_{j_1,\ldots,j_k}\termj.
\end{equation}
In the monomial, if any color appears once, then by independence and
mean zero condition, $\mathbb{T}_{\mathbf j}=0$ for every $\mathbf j$.
Hence, $\widehat{\mu}_n(q)=0$.

So henceforth, assume that each color appearing in the monomial,
appears at least twice.
Now again, if $\mathbf j$ belongs to a circuit which is not color
matched, then $\mathbb{T}_{\mathbf j}=0$.\vadjust{\goodbreak}

Now form the following matrix $M$:
\[
M(L(i,j))=|X_{i_1}(L(i,j))|+|X_{i_2}(L(i,j))|+\cdots+|X_{i_k}(L(i,j))|.
\]
Observe that
\[
|\termj|\leq\E[M(L(j_1,j_2)\cdots M(L(j_k,j_1)].
\]
From Lemma 1 of~\cite{bosesen}, it is known that
the total contribution of all circuits which have at least one three
match, is zero in the limit.

As a consequence of the above discussion, if $k$ is odd, then $\widehat
{\mu}_n(q)\rightarrow0$.
So assume $k$ is even. In that case, we need to consider only circuits
which are pair matched. Further this pair matching must occur within
the same color.
If $\mathbf j$ belongs to any such circuit, then by independence, mean
zero and variance one condition,
$\mathbb{T}_{\mathbf j}=1$.

Then using all the facts established so far,
\begin{eqnarray*}\lim_{n\rightarrow\infty}\widehat{\mu}
_n(q)&=&\lim
_{n\rightarrow\infty}\frac{1}{n^{k/2+1}}\mathop{\sum_{\pi\dvtx
\pi\ \mathrm{pair}\ \mathrm{matched}}}_{\mathrm{within}\ \mathrm{colors}}\E\bigl[ X_{i_1}({L(\pi(0), \pi(1))})
\cdots\\
&&\hspace*{131.2pt}{}\times
X_{i_k}\bigl({L\bigl(\pi(k-1), \pi(k)\bigr)}\bigr)\bigr]\\
&=&\lim_{n\rightarrow\infty}\frac{1}{n^{k/2+1}}\sum_{w \in
\mathit{CW}_k(2)}\sum
_{\pi\in\Pi_{C_q}(w)}\E\bigl[ X_{i_1}({L(\pi(0), \pi(1))}) \cdots\\
\\
&&\hspace*{147.6pt}{}\times
X_{i_k}\bigl({L\bigl(\pi(k-1), \pi(k)\bigr)}\bigr)\bigr]\\
&=&\sum_{w \in \mathit{CW}_k(2)}p_{C_q}(w).
\end{eqnarray*}
The last claim in part (a) follows since
\[
\sum_{w \in \mathit{CW}_{2k}(2)}p_{C_q}(w)= \sum_{w \in
\mathit{CW}_{2k}(2)}p(\psi_q(w))\leq\sum_{w \in W_{2k}(2)}p(w) \leq\frac
{(2k)!\Delta(L)^k}{k!2^k}.
\]
The last inequality above is shown in~\cite{bosesen}.

(b) For part (b), the following notions will be useful: $l$ circuits
$\pi
_1,\pi_2,\ldots,\pi_l$ are said to be \textit{jointly matched} if each
$L$-value occurs at least twice across all circuits. They are said to
be \textit{cross matched} if each circuit has at least one $L$-value
which occurs in at least one of the other circuits.
We can write
%
%
\begin{equation}\label{fourthone}\E[|\widetilde{\mu}_n(q)-\widehat
{\mu}
_n(q)|^4]=\frac{1}{n^{2k+4}}\sum_{\pi_1,\pi_2,\pi_3,\pi_4}\E\Biggl[
\prod
_{l=1}^4 (\mathbb X_{\pi_l}-\E\mathbb X_{\pi_l})\Biggr],
\end{equation}
where
\[
\mathbb X_{\pi} = X_{i_1}({L(\pi(0), \pi(1))}) \cdots
X_{i_k}\bigl({L\bigl(\pi(k-1), \pi(k)\bigr)}\bigr).
\]
If $(\pi_1,\pi_2,\pi_3,\pi_4)$ are not jointly matched, then one of the
circuits, say $\pi_j$, has an $L$ value which does not occur anywhere
else. Also note that $\E\mathbb X_{\pi_j}=0$. Hence, using independence
%
%
\begin{equation}\label{fourthtwo}\E\Biggl[ \prod_{l=1}^4 (\mathbb X_{\pi
_l}-\E\mathbb X_{\pi_l})\Biggr]=\E\Biggl[ \mathbb X_{\pi_j}\prod_{l=1,l\neq j}^4
(\mathbb X_{\pi_l}-\E\mathbb X_{\pi_l})\Biggr]=0.
\end{equation}
If $(\pi_1,\pi_2,\pi_3,\pi_4)$ is jointly matched but is not cross
matched then one of the circuits, say $\pi_j$ is only self-matched,
that is, none of the $L$-values is shared by the other circuits. Then
by independence,
%
%
\begin{equation}\label{fourththree}\quad
\E\Biggl[ \prod_{l=1}^4 (\mathbb X_{\pi
_l}-\E\mathbb X_{\pi_l})\Biggr]=\E\Biggl[ (\mathbb X_{\pi_j}-\E\mathbb X_{\pi
_j})\prod_{l=1,l\neq j}^4 (\mathbb X_{\pi_l}-\E\mathbb X_{\pi_l})\Biggr]=0.
\end{equation}
Since $\{X_{i,n}\}_{1\leq i\leq n}$ satisfy Assumption~\ref{AssumpI},
$ \E[ \prod
_{l=1}^4 (\mathbb X_{\pi_l}-\E\mathbb X_{\pi_l})]$ is uniformly bounded
over all $(\pi_1,\pi_2,\pi_3,\pi_4)$.

The arguments given in~\cite{bry} for Toeplitz and Hankel matrices can
be extended to our set up easily to yield the following: let $Q_{k,4}$
be the number of quadruples of circuits $(\pi_1,\pi_2,\pi_3,\pi_4)$ of
length $k$ such that they are jointly matched and cross matched with
respect to $L$. If $L$ satisfy Property~\ref{properB}, then there
exists a constant
$K$ such that
$Q_{k,4}\leq Kn^{2k+2}$.
Using this, and~\mbox{(\ref{fourthone})--(\ref{fourththree})},
\[
\E[|\widetilde{\mu}_n(q)-\widehat{\mu}_n(q)|^4]\leq K\frac
{n^{2k+2}}{n^{2k+4}}=O(n^{-2}).
\]
Now by an easy application of Borel--Cantelli lemma $\widetilde{\mu}
_n(q)$ converges almost surely.

\section*{Acknowledgments}
We thank the anonymous referees for their constructive comments and
valuable suggestions. That the reverse circulant matrices are half
commuting as well as some important references were pointed out by the
referees. We are grateful to Roland Speicher for his comments and
suggestions.

%

%
\printaddresses

\end{document}